\documentstyle[12pt,amscd]{amsart}
\newtheorem{thm}{Theorem}
\newtheorem{lem}{Lemma}
\newtheorem{prop}{Proposition}
\newtheorem{cor}{Corollary}
\newtheorem{defn}{Definition}
\theoremstyle{definition}

\def\nat{{\Bbb N}}
\def\A{{\cal A}}
\def\B{{\cal B}}

\def\M{{\cal M}}

\def\chix{\raise.5ex\hbox{$\chi$}}

\begin{document}
\baselineskip=18pt
\title[Absolutely Summing Operators]{ Absolutely Summing Operators on 
non commutative $C^*$-algebras and applications}
\author{Narcisse Randrianantoanina}
\address{Department of Mathematics, The University of Texas at Austin,
Austin, TX 78712-1082}
\email{nrandri@@math.utexas.edu}


\maketitle

\begin{abstract}
Let $E$ be a Banach space that does not contain any copy of $\ell^1$ and $\A$
be a non commutative $C^*$-algebra. 
We prove that every absolutely summing operator from $\A$ into $E^*$ is
compact, thus answering a question of Pe\l czynski.

As application, we show that if $G$ is a compact metrizable abelian group
and $\Lambda$ is a Riesz subset of its dual then every countably additive
 $\A^*$-valued
measure with bounded variation and whose Fourier transform is supported by
$\Lambda$ has relatively compact range. Extensions of the same result to 
symmetric spaces of measurable operators are also
 presented.
\end{abstract}

\section{Introduction}

It is a well known result that every absolutely summing operator
from a $C(K)$-space into a separable dual space is compact.
More generally if $F$ is a Banach space with the complete continuity property
 (CCP) then  every absolutely summing operator from any $C(K)$-spaces
  into $F$ is
compact (see \cite{EM5}).

It is the intention of the present note to study extensions of the above
results in the setting of $C^*$-algebras, i.e., replacing the $C(K)$-spaces
above  by a
general non commutative $C^*$- algebra. Typical examples of Banach
spaces with the CCP are dual spaces whose preduals do not contain $\ell^1$.
Our main result is  that  if $E$
is a Banach space that does not contain any copy of
 $\ell^1$ and $\A$ is a 
$C^*$-algebra then every absolutely summing operator from $\A$ into $E^*$
 is compact. This answered positively  the following
 question  raised  by  Pe\l czynski (see \cite{HN} Problem 3. P.20): 
Is every absolutely summing operator from a 
non commutative $C^*$-algebra into a Hilbert space
 compact?
 This result is also used to  study  relative compactness of
  range of  countably additive vector
 measures with values in  duals of  non commutative
$C^*$-algebras.
In \cite{EDG}, Edgar introduced new types of Radon-Nikodym properties 
associated with Riesz subsets of countable discrete group (see the definition 
below) as generalization of the usual Radon-Nikodym property (RNP) and the 
Analytic Radon-Nikodym property (ARNP). These properties were extensively
studied in \cite{DO2} and \cite{DO4}. In \cite{DO2}, it was shown that if 
$\Lambda$ is a Riesz subset of a countable discrete group then $L^1[0,1]$ has
the type II-$\Lambda$-RNP. In the other hand, Haagerup and Pisier
showed in \cite{HP} that non commutative $L^1$-spaces have the ARNP so it is
 a natural question to ask if non commutative $L^1$-spaces have
 the type II-$\Lambda$-RNP for any Riesz subset. In this direction, 
  we obtain (as a consequence of our main result) that if a 
countably additive vector
measure of bounded variation is defined on the $\sigma$-field of
 Borel subsets of a compact
metrizable abelian group and
 takes its values in a dual of a $C^*$-algebra
 then its range is relatively compact provided that its Fourier
transform is supported by a Riesz subset of the dual group.

Our terminology and notation are standard.
We refer to \cite{D1} and \cite{DU} for definitions from Banach space
theory and \cite{TAK}, \cite{N} and \cite{FK} for basic properties from the
theory of operator algebras and non-commutative integrations.

\section{Preliminary Facts and Notations}

We recall some definitions and well known facts which we use in the sequel.

Let $\A$ be a $C^*$-algebra, we denote by $\A_h$ the set of Hermitian
(self adjoint) elements of $\A$.

\begin{defn}
Let $E$ and $F$ be Banach spaces and $0<p<\infty$.
An operator $T:E\to F$ is said to be absolutely
$p$-summing (or simply $p$-summing) if there exists $C$
such that
$$\biggl( \sum_{i=1}^\ell \|Te_i\|^p\biggr)^{1/p} \le
C\max \Bigl\{ \sum |\langle e_i,e^*\rangle |^p, \|e^*\|\le 1\Bigr\}^{1/p}$$
\end{defn}

The following class of operators was introduced by 
Pisier in \cite{PIS3}
as extension of the $q$-summing operators in the setting of $C^*$-algebras.

\begin{defn}
Let $\A$ be a $C^*$-algebra and $F$ be a Banach space, $0<q<\infty$.
An operator $T:\A\to F$ is said to be $q$-$C^*$-summing if there exists a
constant $C$ such that for any finite sequence $(A_1,\ldots,A_n)$ of
Hermitian elements of $\A$ one has
$$\biggl( \sum_{i=1}^n \|T(A_i)\|^q\biggr)^{1/q}
\le C\Big\| \Big( \sum_{i=1}^n |A_i|^q\Bigr)^{1/q}\Big\|_{\A}\ .$$
\end{defn}
The smallest constant $C$ for which the above inequality holds is denoted
by $C_q(T)$.
It should be noted that if the $C^*$-algebra $\A$ is commutative then
every $q$-$C^*$-summing operator from $\A$ into any Banach space is
 $q$-summing.
The following extension of the classical
 Pietsch's factorization  theorem (\cite{D1})
was obtained by Pisier (see Proposition~{1.1} of \cite{PIS3}).

\begin{prop}
If $T:\A\to F$ is  a $q$-$C^*$-summing operator
then there exists a positive linear
form $f$ of norm less than $1$ such that
$$\|Tx\| \le C_q (T) \{f(|x|^q)\}^{1/q}\ ,\ \text{for every}\ x\in \A_h.$$
\end{prop}

Let $\M$ be a von-Neumann algebra and $\M_*$ be its predual. We recall that
a functional $f$ on $\M$ is called normal if it belongs to $\M_*$. 
In \cite{PIS4}, it was shown 
that for the case of von-Neumann algebra
and the operator $T$ being weak* to weakly continuous then
 the positive linear form on the above proposition can be chosen
to be normal; namely we  have the following lemma (see Lemma~4.1 of
 \cite{PIS4}).

\begin{lem}
Let $T:\M\to F$ be a $1$-$C^*$-summing operator.
 If $T$ is weak* to  weakly continuous
 then there exists a linear form
$f\in \M_*$ with $\|f\|\le1$ such that
\begin{equation}
\|Tx\| \le C_1 (T) f(|x|)\ ,\ \text{for every}\ x\in \M_h.
\tag 1
\end{equation}
\end{lem}

 For the next lemma, we recall that for $x \in \M$ and $f \in \M^*$,
$xf$ (resp. $fx$) denotes the element of $\M^*$ defined by 
 $xf(y)=f(yx)$ (resp. $fx(y)=f(xy)$) for all $y \in \M$.

\begin{lem}
Let $f$ be a positive linear form on $\M$. For 
every $x \in \M$,
\begin{equation}
f\left( \left(\frac{xx^* +x^*x}{2}\right)^{1/2} \right) \leq 2 \| xf +fx\|_{\M^*}.
\tag 2
\end{equation}
\end{lem}
 
\begin{pf}
Assume first that $x \in \M_h$. In this case
 $(xx^* + x^*x)/2=|x|^2$. The operator $x$ can be decomposed
as $x= x_{+} - x_{-}$ where $x_{+},x_{-} \in \M^{+}$ and $x_{+}x_{-}=0$.
There exists a projection $p \in \M$ such that $px_{-}=x_{-}p=x_{-}$ and
$(1-p)x_{+}=x_{+}(1-p)=x_{+}$. This yields the following estimates:
\begin{align*}
f(|x|)&= f(x_{+} + x_{-})=f(x_{+}) +f(x_{-}) \cr
 &=\frac{1}{2}(xf +fx)(1-p) + \frac{1}{2}(xf + fx)(p) \cr
 &\leq \frac{1}{2} \left(\|xf +fx\| \|1-p\| + \| xf +fx\| \|p\|\right) \cr
 &\leq \|xf +fx\|.
\end {align*}
For the general case, fix $x \in \M$. Let $a=(x+x^*)/2$ and 
$b=(x-x^*)/{2i}$. Clearly $x= a +ib$ and 
$\left( (xx^* +x^*x)/2\right)^{1/2}= |a| +|b|$. Using the Hermitian case,
we get:
\begin{align*}
f(|a| +|b|) &\leq \| af +fa\| + \|bf +fb\| \cr
            &\leq \|xf +fx\| + \|x^*f +fx^*\|;
\end{align*}
but since $f\geq 0$,
 \begin{align*}
 \|x^*f +fx^*\| &=\sup\{ |f(sx^* +x^*s)|;\ s \in \M, \|s\|\leq 1\} \cr
                &=\sup\{ |f^*(xs^* + s^*x)|;\ s \in \M, \|s\| \leq 1\} \cr
                &=\sup\{|f(xs^* + s^*x)|;\ s \in \M, \|s\| \leq 1\}\cr
                &=\|xf +fx\|,
\end{align*}
which completes the proof of the lemma. 
\end{pf}
\section{Main Theorem}

\begin{thm}
Let $\A$ be a $C^*$-algebra,
 $E$ be a Banach space that does not contain any copy of $\ell^1$
and $T:\A\to E^* = F$ be a
$1$-summing  operator then $T$ is compact.
\end{thm}

We will divide the proof into two steps. First we will assume that the
 $C^*$-Algebra $\A$ is a $\sigma$-finite von-Neumann algebra and the operator 
$T$ is weak* to weakly continuous; then we will show that the general can be
 reduced to this case. We refer to \cite{TAK} P.~78 for the definition of
$\sigma$-finite von-Neumann algebra. 

\begin{prop} Let $\M$ be a $\sigma$-finite von-Neumann algebra.
 Let $T: \M \to E^*$ be a  weak* to weakly continuous 1-summing operator
then $T$ is compact.
\end{prop}

\begin{pf} The operator $T$ being weak* to weakly continuous 
and 1-summing, there exist
a constant $C=C_1(T)$ and a normal positive functional $f$ on $\M$ such that
 $$\|Tx\| \leq C f(|x|)\, \ \text{for every }\ x \in \M_h.$$
Since the von-Neumann algebra $\M$ is $\sigma$-finite, there exists a faithful
normal state $f_0$ in $\M_*$ (see \cite{TAK} Proposition~II-3.19).
 Replacing $f$ by $f+f_0$, we can assume that the
functional $f$ on the inequality  above is a faithful normal state and using
Lemma~2, we get
 \begin{equation}
\|Tx\| \leq 2C \|xf + fx\|_{\M^*}\ \text{for every}\ x\in \M.
\tag 3
\end{equation}
We may equip $\M$ with the scalar product by setting for every
 $x,\ y\ \in \M$,
 $$ \langle x,y\rangle =f\left(\frac{xy^* +y^*x}{2}\right).$$
Since $f$ is faithul, $\M$ with $\langle.,.\rangle$ is a pre-Hilbertian.
 we denote  the completion of this space by $L^2(\M,f)$ (or simply $L^2(f)$).

By construction, the inclusion map $J: \M \to L^2(\M,f)$ is bounded
and is one to one ($f$ is faithful).
On the dense subspace $J(\M)$ of $L^2(f)$, we define a map
$\theta:\ J(\M) \to L^2(f)^*$ by $ \theta(Jx)=\langle.,J(x^*)\rangle$.
The map  $\theta$ is clearly linear and is an isometry; indeed for every 
$x \in \M$, $\|\theta(Jx)\|^2=\sup\limits_{\|u\|\le 1}\langle u,J(x^*)\rangle^2
=\langle J(x^*),J(x^*)\rangle=f(x^*x +xx^*)=\|Jx\|^2$.
 So it can be extended to
a bounded map (that we will denote also by $ \theta$) from $L^2(f)$ onto
$L^2(f)^*$.

Let $S= J^* \circ\theta\circ J$. The operator is defined from $\M$ into $\M_*$
 and we  claim that for every $x \in \M$,
 $Sx= xf +fx$. In fact for every $x,\ y \in \M$, we have:
 \begin{align*}
 Sx(y) &=J^*\circ \theta\circ Jx (y)\\
 &=\theta\circ Jx (Jy) \\
 &=\langle J(y),J(x^*) \rangle \\
 &=f(xy +yx)=(xf +fx)(y).
\end{align*}
Notice also that since $f$ is normal, the functionals 
$xf$ and $fx$ are both normal for every $x \in \M$; therefore
$S(\M)\subset \M_*$. Also since $J$ is one to one, $J^*$ has weak* dense 
range. The latter with the facts that both $J$ and $\theta$ have dense ranges
imply that $S(\M)$ is weak* dense in $\M^*$ so $S(\M)$ is (norm) dense in
$\M_*$.

Let us now define a map $L:\ S(\M) \to E^*$ by
$L(xf +fx)= Tx$\ for every $x \in \M$. The map $L$ is clearly linear and 
one can deduce from (3) that $L$ is bounded  so it can be extended as a
bounded operator (that we will denote also by $L$) from $\M_*$ into 
$E^*$. The above means that $T$ can be factored as follows
$$\matrix
\M &\buildrel{T}\over\longrightarrow & E^* \cr
 \qquad\llap{$\scriptstyle S$}
\searrow&&\nearrow\rlap{$\scriptstyle L$} \qquad \cr
& \M_*
\endmatrix
$$
Taking the adjoints we get
$$\matrix
E &\buildrel{T^*}\over\longrightarrow &\M_* \cr
\qquad\llap{$\scriptstyle L^*$}
\searrow&&\nearrow\rlap{$\scriptstyle S^*$} \qquad \cr
&\M
\endmatrix
$$

To conclude the proof of the proposition, let $(e_n)_n$ be a bounded
sequence in $E$. Since $E \not\hookleftarrow \ell^1$,
we will assume (by taking a subsequence if necessary) that $(e_n)_n$
is weakly Cauchy.
We will show that $(T^*(e_n))_n$ is norm-convergent.
For that it is enough to prove that if $(e_n)_n$ is a weakly null sequence
in $E$ then $(\|T^* e_n\|)_n$ converges to zero.

Let $(e_n)_n$ be a weakly null sequence in $E$,
$(L^*(e_n))_n$ is a weakly null sequence in $\M$.
This implies that  $((L^* (e_n))^*)_{n\ge1}$ (the sequence of the
adjoints of the $L^*(e_n)$'s) is weakly
null in $\M$.

Since $T$ is 1-summing, it is a Dunford-Pettis operator
(i.e takes weakly convergent sequence into norm-convergent 
sequence). Hence
$$\lim_{n\to\infty} \|T( (L^* e_n)^* )\|_{E^*} = 0\ .$$
In particular, since $(e_n)_n$ is a bounded sequence in $E$, we have
$$\lim_{n\to\infty}\langle T( (L^* e_n)^* ),e_n\rangle =0$$
but
\begin{align*}
\langle T( (L^* e_n)^* ),e_n\rangle
& = \langle LS((L^* e_n)^* ),e_n\rangle \cr
& = \langle  S((L^* e_n)^* ), L^*e  _n\rangle \cr
& = \langle \theta\circ J((L^* e_n)^*), J(L^* e_n)\rangle\cr
& = \langle J(L^* e_n), J(L^* e_n) \rangle_{L^2(f)} \cr
& = \|J(L^* e_n)\|_{L^2(f)}.
\end{align*}
So $\|J(L^* e_n)\|_{L^2(f)}\to 0$ as $n\to\infty$ and therefore
since $T^*=S^*\circ L^* = J^*\circ\theta \circ J\circ L^*$,
we get that 
  $\lim_{n\to\infty} \|T^* e_n\|=0$.

This shows that $\overline{T^* (B_E)}$ is compact and since $B_E$ is weak*
dense in  $B_{E^{**}}$ and $T^*$ is weak* to weakly continuous,
$T^* (B_{E^{**}}) \subseteq \overline{T^* (B_E)}$
so  $T^*$ (and hence $T$) is compact.
The proposition is proved.
\end{pf}

To complete the proof of the theorem, let $\A$ be a $C^*$-algebra and
$T: \A \to E^*$ be a 1-summing operator. The double dual $\A^{**}$ of $\A$
is a von-Neumann and $T^{**}:\ \A^{**} \to E^*$ is 1-summing. Let $(a_n)_n$ be
a bounded sequence in $\A^{**}$. If we denote by $\M$  the von-Neumann algebra
generated by $(a_n)_n$ then the predual $\M_*$ of $\M$ is separable and 
therefore the von-Neumann algebra $\M$ is $\sigma$-finite. Moreover, if we
 set $I: \M \to \A^{**}$ the inclusion map then $I$ is weak* to weak*
 continuous. Hence $\M$ and $T^{**} \circ I$ satisfy the conditions of 
Proposition~2 so $T^{**}\circ I$ is compact and since the sequence 
$(a_n)_n$ is arbitrary, the operator $T^{**}$ (and hence $T$) is compact.
\qed

\vskip .3truein
\noindent
{\bf Remark.} It should be noted that for the  proof of Proposition~2, 
we only require the
operator $T$ to be $C^*$-summing and Dunford-Pettis so
the conclusion of  Proposition~2 is still
valid for $C^*$-summing operators that are Dunford-Pettis.

\section{Applications to vector measures}

In this section we will provide some  applications of the main theorem 
to study
range of countably additive vector measures with values in duals of 
$C^*$-Algebras.

The letter $G$ will denote a compact metrizable abelian group,
$\widehat{G}$ its dual, $\cal B(G)$ is
the $\sigma$-algebra of the Borel subsets of $G$, and $\lambda$ the
normalized Haar measure on $G$.

Let $X$ be a Banach space and $1\leq p\leq \infty$, we will denote by 
$L^p(G,X)$ the usual Bochner spaces for the measure space
 $(G,\cal B(G),\lambda)$; $M(G,X)$ the space of $X$-valued countably additive
Borel measures of bounded variation; $C(G,X)$ the space of $X$-valued
 continuous functions  and $M^\infty(G,X)=\{\mu \in M(G,X),\ |\mu|\le
C\lambda\ \text{for some}\  C>0\}$.

If $\mu \in M(G,X)$,
we recall that the Fourier transform  of $\mu$ is a map $\hat\mu$ from
$\widehat G$ into $X$ defined by
$\widehat \mu (\gamma) = \int_G \bar{\gamma} \ d\mu\ \text{ for }\ 
\gamma\in\widehat G$.

For $\Lambda \subset  \widehat G$, we will use the following notation:
\begin{align*}
 L_{\Lambda}^p(G,X) &=\{f \in L^p(G,X),\ \hat{f}(\gamma)=0 \ \text{for all}\
\gamma \notin \Lambda\} \cr
 C_{\Lambda}(G,X)&=\{f \in C(G,X),\ \hat{f}(\gamma)=0 \ \text{for all}\
\gamma \notin \Lambda\} \cr
 M_{\Lambda}(G,X)&=\{\mu \in M(G,X),\ \hat\mu(\gamma)=0 \ \text{for all}\
\gamma \notin \Lambda\} \cr
 M_{\Lambda}^\infty(G,X)&=\{\mu \in M^\infty(G,X),\ \hat\mu(\gamma)=0\ 
\text{for all}
\ \gamma \notin \Lambda\}.
\end{align*} 
We also recall that $\Lambda \subset \widehat G$ is called a Riesz subset if
$M_\Lambda(G)=L_\Lambda^1(G)$.
We refer to \cite{RU} and \cite{GO}
 for detailed discussions  and examples of Riesz subsets of dual groups.

The following Banach space properties were introduced by Edgar \cite{EDG}, and
Dowling \cite{DO2}.
\begin{defn} Let $\Lambda$ be a Riesz subset of $\widehat G$. 
A Banach space $X$ is said to have type I-$\Lambda$-Radon Nikodym Property
(resp. type II-$\Lambda$-Radon Nikodym property) if
 $M_{\Lambda}^\infty(G,X)= L_{\Lambda}^\infty(G,X)$
 (resp. $M_{\Lambda}(G,X)=L_{\Lambda}^1(G,X)$).
\end{defn}

Our next result deals with property of dual of $C^*$-algebras related to the
 types of Radon-Nikodym properties defined above. 
\begin{thm}
Let  $\Lambda$ be a Riesz
subset of $\widehat G$ and $\A$ be a $C^*$-Algebra.
If $F:\B(G)\to \A^*$ is  a countably additive  measure
 with bounded variation that satisfies
$\widehat F(\gamma)=0$ for $\gamma\notin \Lambda$ then
the range of $F$ is a relatively compact subset of $\A^*$.
\end{thm}

\begin{pf}
Let $F:\B(G) \to \A^*$ be a measure with bounded variation and
$\widehat F(\gamma)=0$ for $\gamma\notin \Lambda$.
Let $S:C(G) \to \A^*$ be the operator defined by $Sf=\int f\,dF$.
Since $F$ is of bounded variation, the operator $S$ is integral (see
\cite{DU} Theorem~IV-3.3 and Theorem~IV-3.12) 
and therefore $S^*: \A^{**}\to (C(G))^*$ is also integral.
Now since $\widehat F(\gamma)=0$ for $\gamma\notin \Lambda$, if we
denote by $\Lambda' = \{\gamma\in\widehat G,\bar\gamma\notin \Lambda\}$
then $S(\gamma)=0$ for all $\gamma\in \Lambda'$ and therefore we have
the following factorization
$$\matrix
C(G)\quad \buildrel S\over\longrightarrow&\A^*\cr
\llap{$\scriptstyle q$}\big\downarrow\qquad
&\nearrow\!\!{\scriptstyle L}\qquad\qquad\cr
C(G)/C_{\Lambda'}(G)&
\endmatrix$$
where $q$ is the natural quotient map.
Taking the adjoints, we get
$$\matrix
\A^{**}&\buildrel {S^*}\over\longrightarrow&(C(G))^*\cr
\qquad\llap{$\scriptstyle L^*$}\searrow&&\nearrow\rlap{$\scriptstyle q^*$}
\qquad\cr
&M_\Lambda (G)&
\endmatrix$$
Since $q^*$ is the formal inclusion and $S^*$ is 1-summing, the operator
$L^*$ is 1-summing.
The assumption $\Lambda$ being  a Riesz subset implies that
 $M_\Lambda(G)=L_\Lambda^1(G)$ is
a separable dual (in particular its predual does not contain $\ell^1$).
So by Theorem~1, $L^*$ (and hence $S$) is compact.
This proves that the range of the representing measure $F$ of $S$ 
 is relatively compact (see \cite{DU} Theorem~II-2.18).
\end{pf}

Our next result is a generalization of Theorem~2 for the case of symmetric
spaces of measurable operators.

Let  $(\M,\tau)$ be  a semifinite von-Neumann algebra
acting on a Hilbert space $H$. Let  $\tau$ be  a
distinguished faithful normal semifinite trace on $\M$.

Let $\overline{\M}$ be the space of all measurable operators with respect
to $(\M,\tau)$ in the sense of \cite{N};
for $a\in \overline{\M}$ and $t>0$, the $t^{th}$-s-number (singular number)
of $a$ is defined by
$$\mu_t (a) = \inf \{\|ae\| :e\in \M\text{ projection with }
\tau(I-e)\le t\}\ .$$
The function $t\mapsto \mu_t (a)$ defined  on $(0,\tau(I))$ will be denoted
 by
$\mu(a)$. This is a positive non-increasing function on $(0,\tau(I))$.
We refer to \cite{FK} for complete detailed study  of $\mu(a)$.

Let $E$ be a rearrangement invariant Banach function space on $(0, \tau(I))$
(in the sense  of \cite{LT}).
We define  the symmetric space$E(\M,\tau)$ of measurable 
operators  by setting
$$E(\M,\tau) = \{ a\in \overline{\M} ;\  \mu (a) \in E\}$$
and $\|a\|_{E(\M,\tau)} = \|\mu(a)\|_E$.

It is well known that $E(\M,\tau)$ is a Banach space and if
 $E= L^p (0,\tau(I))$ $(1\le p\le \infty)$ then $E(\M,\tau)$ coincide with
the usual non-commutative $L^p$-spaces associated with the von-Neumann
algebra $\M$.
 The space $E(\M,\tau)$ is often referred  as the non-commutative version 
of the function space $E$. Some  Banch space properties of these spaces 
 can be found in \cite{CS}, \cite{DDP1} and \cite{X}.

For the case where the trace $\tau$ is finite, we obtain the 
following generalization  of Theorem~2 for symmetric  spaces of measurable
operators.
 
\begin{cor}
Assume that $\tau$ is finite.
Let $E$ be a rearrangement invariant
 function space on $(0,\tau(I))$ that does not contain 
$c_0$ and $\Lambda$ be a Riesz subset of $\widehat G$.
  Let  $F: \B(G)\to  E(\M,\tau)$  be  a countably additive measure
with bounded variation and $\widehat{F}(\gamma)=0$  for every
 $\gamma \notin \Lambda$ then the range of $F$ is relatively compact.
\end{cor}

\begin{pf}
We will begin by reducing the general case to the case where $E(\M,\tau)$
 is separable.
Since $\B(G)$ is countably generated, the  range of $F$ is separable.
Choose $(A_n)_n\subset \B(G)$ so that $\{F(A_n),n\ge1\}$ is dense in
$\{F(A),A\in\B(G)\}$.
Let $\widetilde \M$ be the von-Neumann algebra generated $I$ and $F(A_n)$
($n\ge1$) and $\tilde\tau$ the restriction of $\tau$ in $\widetilde \M$.
Clearly $E(\widetilde \M,\tilde \tau)$ is a closed subspace of $E(\M,\tau)$
and $F(A) \in E(\widetilde\M, \tilde \tau)$ for all $A\in \B(G)$.
Moreover   
 the space $E(\widetilde \M,\tilde \tau)$ is separable (see Lemma 5.6 of
\cite{X}). So without loss of generalities we will assume that $E(\M,\tau)$
is separable.
 It is a well known fact that $E(\M,\tau)$ is contained in
 $L^1(\M,\tau) + \M$ and since $\tau$ is finite,
 $E(\M,\tau)\subset L^1(\M,\tau)$.
  Let
 $J: E(\M,\tau) \to L^1(\M,\tau)$ be the formal
inclusion. The measure $J\circ F$ is of bounded variation and 
$\widehat{J\circ F}(\gamma)=J(\widehat F(\gamma))\ \text{for every}\
 \gamma \in \widehat G$. One can conclude from Theorem~2 that the
range of $J\circ F$ is relatively compact in $L^1(\M,\tau) $.

To show that the range of $F$ is relatively compact, fix
$h: G\to E(\M,\tau)^{**}$ a weak*-density of $F$ with respect
to the Haar measure $\lambda$ (see \cite{DIU}). We have for each $A \in \B(G)$,
$$F(A)=\text{weak*}-\int_A h(t)\ d\lambda(t)$$
and $$|F|(A)= \int_A \|h(t)\|\ d\lambda(t).$$
For each $N\in\nat$, let $A_N=\{ t \in G, \|h(t)\| \leq N\}$ and $F_N$ the
 measure defined by $F_N (A) = F(A \cap A_N)$ for all $A\in \B(G)$. Clearly
$|F_N| \leq N\lambda \ \text{for every}\  N \in \nat$.

Define $T_N: L^1(G) \to E(\M,\tau)$ by 
$T_N(f)=\int f(t)\ dF_N(t)\ \ \text{for every}\   f \in L^1(G)$.
The operator $T_N$ is bounded and
 we  claim that $T_N$ is
 Dunford-Pettis; for that  notice that since  the range of $J\circ F$ is
relatively  compact  so is the range of $J\circ F_N$ and therefore
the operator $J\circ T_N$ is a Dunford-Pettis operator. The space
 $E(\M,\tau)$  is separable 
 and $J$ is a semi-embedding (see Lemma 5.7 of \cite{X}) so $J$  is
a $G_\delta$-embedding (see \cite{BR} Proposition~1.8)
 and one can deduce from
Theorem~II.6 of \cite{GR}, that $T_N$ is a Dunford-Pettis operator.
 Hence the range of $F_N$ is relatively compact.
Now since
 $$\lim_{N\to \infty}\|F-F_N\|=\lim_{N\to \infty}\int_{G\setminus A_N}\|h(t)\|\
 d\lambda(t)=0,$$
the range of $F$ is relatively compact. 
\end{pf}

\vskip .3truein
Let us finish by asking the following question:

\noindent{\bf Question:} Do non-commutative $L^1$-spaces have
 type II-$\Lambda$-RNP for any Riesz set $\Lambda$?

In light of Theorem~2, the result of Haagerup and Pisier (\cite{HP}) and so
many properties that have been generalized from classical $L^1$-spaces to
non-commutative $L^1$-spaces, one tends to conjucture that the answer of
 the above question is affirmative. 

\vskip .3truein
\noindent{\bf Ackowlegements:} I  would like to thank Professor 
 G. Pisier for many valuable suggestions conserning this work.


\end{document}